\documentclass[pamm,a4paper,fleqn]{w-art}
\usepackage{times,cite,w-thm}
\usepackage[T1]{fontenc}
\usepackage[utf8]{inputenc}
\numberwithin{equation}{section}
\usepackage[ngerman,english]{babel}
\usepackage{graphicx}
\usepackage{booktabs}
\usepackage{caption}
\usepackage{subcaption}

\definecolor{royalpurple}{rgb}{0.47, 0.32, 0.66}

\usepackage{ifthen}%
\newcommand{\mycomment}[1]{%
  \ifthenelse{\isodd{\value{page}}}{%
    \normalmarginpar%
    \marginpar{\tiny {#1}}%
  }{%
    \reversemarginpar%
    \marginpar{\tiny {#1}}%
  }}%

\newcommand{\takeout}[1]{}

\begin{document}
%% \def\leftmark{Session title}
%%
%%    The information for the title page will be placed between
%%    \begin{document} and \maketitle. The order of most entries
%%    is determined by the class file and cannot be changed by
%%    rearranging them. The maketitle command follows after the
%%    abstract.
%%
%%    The following commands will be updated by the publisher:
%%
%%    \renewcommand{\copyrightyear}{2016}
%%    \DOIsuffix{pamm.20161zzzz}
%%    \Volume{16}
%%    \Year{2016}
%%    \pagespan{1}{}
%%
%%    The short title is optional:

\TitleLanguage[EN]
\title[Monolithic Multi-level Overlapping Schwarz]{Monolithic Multi-level Overlapping Schwarz Solvers for Fluid Problems}

%% Please delete not needed author entries.
%% Information for the first author.
\author{\firstname{Stephan} \lastname{Köhler}\inst{1}%
\footnote{e-mail \ElectronicMail{stephan.koehler@math.tu-freiberg.de},
  phone +49\,3731\,39\,3188}}
\address[\inst{1}]{\CountryCode[DE]INMO, Technische Universität Bergakademie Freiberg, 09596 Freiberg}
\address[\inst{2}]{\CountryCode[DE]University Computing Centre (URZ), Technische Universität Bergakademie Freiberg, 09596 Freiberg}
\address[\inst{3}]{\CountryCode[DE]Center for Efficient High-Temperature Processes and Materials Conversion (ZeHS), Technische Universität Bergakademie Freiberg, 09596 Freiberg}

%%
%%    Information for the second author
\author{\firstname{Oliver} \lastname{Rheinbach}\inst{1,2,3}%
  \footnote{e-mail \ElectronicMail{oliver.rheinbach@math.tu-freiberg.de},
    phone +49\,3731\,39\,3279}}
%% \address[\inst{2}]{\CountryCode[FR]Second address}
%%
%%
%%    \dedicatory{This is a dedicatory.}
%%
%%    Abstract is required.
\AbstractLanguage[EN]
\begin{abstract}
    Additive overlapping Schwarz Methods are iterative methods of the domain decomposition type for the solution of partial differential equations. Numerical and parallel scalability of these methods can be achieved by adding coarse levels. A successful coarse space, inspired by iterative substructuring, is the generalized Dryja–Smith–Widlund (GDSW) space.  In~\cite{heinlein2019monolithic}, based on the GDSW approach, two-level monolithic overlapping Schwarz preconditioners for saddle point problems were introduced.  We present parallel results up to 32\,768 MPI for the solution of incompressible fluid problems for a Poiseuille flow example on the unit cube and a complex extrusion die geometry ranks using a two- and a three-level monolithic overlapping Schwarz preconditioner.  These results are achieved through the combination of the additive overlapping Schwarz solvers implemented in the Fast and Robust Overlapping Schwarz (FROSch) library~\cite{Heinlein:2018:FPI}, which is part of the Trilinos package ShyLU~\cite{Shylu}, and the FEATFLOW library~\cite{featflow-website} using a scalable interface for the efficient coupling of the two libraries. This work is part of the project StroemungsRaum - Novel Exascale-Architectures with Heterogeneous Hardware Components for Computational Fluid Dynamics Simulations, funded by the German Bundesministerium für Forschung, Technologie und Raumfahrt BMFTR (formerly BMBF) as part of the program on New Methods and Technologies for Exascale Computing (SCALEXA).
\end{abstract}
%% maketitle must follow the abstract.
\maketitle                   % Produces the title.

\section{Introduction}
We consider monolithic overlapping Schwarz preconditioners of the generalized Dryja-Smith-Widlund (GDSW) type~\cite{Dohrmann:2008:FEM,Dohrmann:2008:DDL} for the application to fluid problems.
This work is part of the StroemungsRaum project, financed by the BMFTR (formerly BMBF) under the SCALEXA program. SCALEXA is the German initiative to develop software in high-performance computing (HPC) in the exascale computing era.
StroemungsRaum aims to extend the CFD software FEATFLOW~\cite{featflow-website} with highly scalable domain decomposition methods and multigrid solvers tailored for heterogeneous CPU-GPU systems. By leveraging parallel near-hardware implementations, the project targets efficient CFD simulations on future exascale architectures. %Time-parallel and time-simultaneous techniques will further boost parallel scalability.

This paper is concerned with the enhancement of the FEATFLOW software library by scalable domain decomposition methods based on the FROSch~\cite{Heinlein:2018:FPI} framework 
with a focus on GDSW with three and more levels~\cite{Roever:2019}.
FROSch is part of the ShyLU package~\cite{Shylu} within the Trilinos software library~\cite{trilinosrepo}, and provides a flexible infrastructure for implementing advanced overlapping domain decomposition preconditioners.

In the numerical simulation of fluid dynamics problems, the underlying discretized systems typically exhibit a saddle point structure. As the geometric complexity increases and finer mesh resolutions are required to resolve detailed flow features or boundary layers, the resulting linear systems can become extremely large. For such systems, direct solvers become infeasible due to their high computational cost and memory consumption.
This motivates the use of iterative methods, particularly Krylov subspace methods, combined with multilevel preconditioners from multigrid and domain decomposition. 

Domain decomposition methods, such as GDSW in the FROSch framework, offer a promising path to scalable preconditioning by enabling the solution process to be distributed across multiple compute nodes, exploiting both spatial locality and parallelism. In this context, multi-level strategies are essential to maintain scalability on emerging architectures with large core counts.

FEAT3 is a modern C++ finite element software developed at TU Dortmund as a modular successor to earlier versions of the FEAT family. It emphasizes highly scalable geometric multigrid solvers, enabling efficient simulation of large-scale PDE systems on modern HPC architectures. Within the StroemungsRaum project, FEAT3 is used in collaboration with the industrial partner IANUS Simulation, who applies it to deliver simulation-as-a-service solutions.

In~\cite{heinlein2019monolithic} the monolithic two-level overlapping Schwarz preconditioner with a generalized Dryja–Smith–Widlund (GDSW) coarse spaces was introduced and its scalability and robustness was demonstrated, see also~\cite{hochmuth2020parallel}.
It is well known that with an increasing number of subdomains the size of the GDSW coarse problem grows and hence the factorization of the coarse problem will  eventually be infeasible.  To overcome this problem a multi-level version of overlapping Schwarz preconditioners was introduced in~\cite{heinlein2018three}, see also~\cite{roever2022multilevel:phd, heinlein2023multilevel}.

In this work, we extend the multi-level approach to the monolithic overlapping Schwarz preconditioner and show weak scaling resulting up to 32\,768 MPI ranks for a Poiseuille flow problem on the unit square and up to 4\,096 MPI ranks for Poiseuille flow on the geometry of an extrusion die.

The paper is arranged as follows.  In section~\ref{pamm-25-koehler-sec:stokes-problem} we introduce our model Stokes problem and give a brief overview of the monolithic overlapping Schwarz preconditioner and the extension to the monolithic multi-level preconditioner.  In section~\ref{pamm-25-koehler-sec:implementation}, we discuss briefly our implementation and in section~\ref{pamm-25-koehler-sec:numerical-results}, we show our numerical results.

\section{Stokes Problem}\label{pamm-25-koehler-sec:stokes-problem}
For our monolithic multi-level overlapping Schwarz preconditioner, we consider the Stokes equation over a computational domain $\Omega\subset \mathbf{R}^{d}$, where $d=2,3$.  The weak formulation of the Stokes equation is given by:  Find $u\in V_{g}$ and $p\in Q$ such that
\begin{equation*}
\begin{aligned}
  \mu \int_{\Omega}\nabla u\colon\nabla v\,\mathrm{d}x &\,\, - \int_{\Omega}\mathrm{div}v\,p\,\mathrm{d}x &=&  \int_{\Omega}fv\mathrm{d}x & & \forall v\in V_{0}, \\
  -\int_{\Omega}\mathrm{div}u\,q\,\mathrm{d}x          &                                                 &=&\,  0                      & & \forall q\in Q_{0}, \\[1.5ex]
                                                                                                                                           u = g, &  \quad x\in\partial\Omega_{D}, & &&  \frac{\partial u}{\partial n} -pn = 0,  &\quad x\in\partial\Omega_{N},
\end{aligned}
\end{equation*}
where $\Omega_{D}$ denotes the Dirichlet boundary, $\Omega_{N}$ denotes the Neumann boundary, and $n$ the outward pointing normal vector.  Note that in our numerical experiments, we only consider problems with an inflow and outflow boundary conditions and therefore the pressure is uniquely determined.  We could handle a non-unique pressure in the same way as considered in~\cite{heinlein2019monolithic}.  For the discretization, we use the FEATFLOW library~\cite{featflow-website} and Q2 elements for the velocity part $u$ and P1 discontinuous for pressure part $p$.  The discretized system can be written as
\begin{equation}
  \label{eq:pamm-feat-frosch-saddle-point}
  \mathcal{A}:=
  \begin{bmatrix} A & B^{T} \\[0.5ex] B & 0 \end{bmatrix}
  \begin{bmatrix} u \\[0.5ex] p \end{bmatrix}
  = \begin{bmatrix} F \\[0.5ex] 0 \end{bmatrix} =: \mathcal{F},
\end{equation}
where the matrix $A$ corresponds to the velocity and $B^{T}$ is the coupling between the pressure and velocity.

\subsection{Monolithic Overlapping Schwarz Preconditioner}

For overlapping Schwarz methods, the computational domain $\Omega$ is decomposed into $N$ nonoverlapping subdomains $\Omega_{i}$, $i=1,\ldots,N$ and we denote the corresponding subdomain with $k$ layers of overlap with $\Omega_{i}'$.  We denote the local saddle point problem of $\Omega_{i}'$ with $\mathcal{A}_{i}$, i.e., $\mathcal{A}_{i}$ is the restriction of $\mathcal{A}$ onto the overlapping subdomain $\Omega_{i}'$.  The GDSW coarse spaces for \eqref{eq:pamm-feat-frosch-saddle-point}, see~\cite{heinlein2019monolithic}, can be written as
\begin{equation}
  \label{eq:pamm-feat-frosch-gdsw}
  \mathcal{B}_{\mathrm GDSW} = \Phi \mathcal{A}_{0}^{-1}\Phi^{T} + \sum_{i=1}^{N}\mathcal{R}_{i}^{T}\mathcal{A}_{i}^{-1}\mathcal{R}_{i},
\end{equation}
where
\begin{equation}
  \label{eq:pamm-feat-frosch-gdsw-details}
  \begin{aligned}
    \mathcal{A}_{0}:= \Phi^{T}\mathcal{A}\Phi, & \qquad & \mathcal{R}_{i}:= \begin{bmatrix} R_{i, u} & 0 \\[0.5ex]  0 & R_{i, p} \end{bmatrix},
  \end{aligned}
\end{equation}
the columns of $\Phi$ are the coarse basis functions, $R_{i, u}$ denotes the restriction of the velocity onto the $i$-th subdomain, and $R_{i, p}$ the restriction of the pressure~\cite{heinlein2019monolithic,hochmuth2020parallel}.

Note that $A_{0}$ is invertible since we have a unique pressure due to the outflow boundary condition and the local problems $\mathcal{A}_{i}$ are invertible since they correspond to problems with Dirichlet boundary conditions for the velocity and pressure.  The coarse basis functions are the energy minimizing extensions of the restriction of the null space of the operators onto the interface $\Gamma$ and can be constructed the following:  We define the interface $\Gamma$ as $\Gamma=\bigcup_{i=i}^{N}\partial\Omega_{i}\setminus\partial\Omega_{D}$.  Furthermore, we divide the interface into vertices and edges and, in 3D, faces.  By reordering of the degrees of freedom (dofs), we can write $\mathcal{A}$ as
\begin{equation*}
  \mathcal{A} = \begin{bmatrix} \mathcal{A}_{II} & \mathcal{A}_{I\Gamma} \\[0.5ex] \mathcal{A}_{\Gamma I} & \mathcal{A}_{\Gamma\Gamma}  \end{bmatrix},
\end{equation*}
where subscript $\Gamma$ denotes of the dofs of the interface and subscript $I$ denotes the remaining dofs.  The matrix $\Phi$ can now be written as
\begin{equation}
  \label{eq:pamm-feat-frosch-gdsw-coarse-basis}
  \Phi:= \begin{bmatrix} \Phi_{I} \\[0.5ex] \Phi_{\Gamma} \end{bmatrix} =
  \begin{bmatrix}
    -\mathcal{A}_{II}^{-1}\mathcal{A}_{I\Gamma}\Phi_{\Gamma} \\[0.5ex] \Phi_{\Gamma}
  \end{bmatrix},
\end{equation}
where the columns of $\Phi_{\Gamma}$ are the prescribed restriction of the null space onto the interface entities.  In detail:  $\Phi_{\Gamma}$ has a row for each dof in $\Gamma$ and each column represents an interface entity (vertex, edge, face), where the entries of the rows are the values of the null space.  Since we have a null space for the velocity and the pressure variable, we can write $\Phi_{\Gamma}$ as
\begin{equation*}
  \Phi_{\Gamma} = \begin{bmatrix} \Phi_{\Gamma, u_{0}} & 0 \\[0.5ex] 0 & \Phi_{\Gamma, p_{0}}\end{bmatrix}.
\end{equation*}
The complete matrix $\Phi$ can now be written as
\begin{equation*}
  \Phi % = \begin{bmatrix} \Phi_{I} \\[0.5ex] \Phi_{\Gamma} \end{bmatrix}
  =
  \begin{bmatrix}
    \Phi_{I,u_{0}, u} & \Phi_{I,u_{0}, p} \\[0.5ex]
    \Phi_{I,p_{0}, u} & \Phi_{I,p_{0}, p} \\[0.5ex]
    \Phi_{\Gamma, u_{0}} & 0 \\[0.5ex]
    0 & \Phi_{\Gamma, p_{0}}
  \end{bmatrix}.
\end{equation*}
For more details and a discussion on the off diagonal blocks $\Phi_{I,u_{0}, p}, \Phi_{I,p_{0}, u}$, see, e.g.,~\cite{heinlein2019monolithic,hochmuth2020parallel}.

\subsection{Monolithic Multi-Level Overlapping Schwarz Preconditioner}

For an increasing number of subdomains, the coarse problem $\mathcal{A}_{0}$ is getting larger and therefore the direct factorization of $\mathcal{A}_{0}$ will be more time and memory consuming.  We follow the idea of the multi-level overlapping Schwarz preconditioners~\cite{heinlein2018three,roever2022multilevel:phd}:  We apply the  monolithic two-level overlapping Schwarz preconditioner recursively to $A_{0}$ and replace $\mathcal{A}_{0}^{-1}$ by the approximation with the resulting preconditioner.  The monolithic three-level overlapping Schwarz preconditioner can be written as
\begin{equation}
  \label{eq:pamm-feat-frosch-gdsw-three-level}
  \begin{aligned}
    \mathcal{B}_{\mathrm GDSW,3-\mathrm{level}} =\ &  \underbrace{\Phi^{(2)} \left(\underbrace{\Phi^{(3)}\left.\mathcal{A}_{0}^{(3)}\right.^{-1}\left.\Phi^{(3)}\right.^{T}}_{\text{3rd level}} + \underbrace{\sum_{i=1}^{N^{(2)}}\left.\mathcal{R}_{i}^{(2)}\right.^{T}\left.\mathcal{A}_{i}^{(2)}\right.^{-1}\mathcal{R}_{i}^{(2)}}_{\text{2nd level}} \right)\left.\Phi^{(2)}\right.^{T}}_{\text{coarse levels}} \\
    & + \underbrace{\sum_{i=1}^{N^{(1)}}\left.\mathcal{R}_{i}^{(1)}\right.^{T}\left.\mathcal{A}_{i}^{(1)}\right.^{-1}\mathcal{R}_{i}^{(1)}}_{\text{1st level}},
  \end{aligned}
\end{equation}
where superscript $(1)$ denotes all parts of the first level, superscript $(2)$ of the second level, and superscript $(3)$ of the third level.  The number $N^{(1)}$ denotes the number of subdomains of the first level and $N^{(2)}$ the number of subdomains on the second level (subregions).  This generalization of the monolithic two-level overlapping Schwarz preconditioner is straight forward and can easily be extended to a multi-level approach.

\section{Implementation}\label{pamm-25-koehler-sec:implementation}
Our implementation of a monolithic multi-level overlapping Schwarz preconditioner is based on the ShyLU/FROSch package~\cite{heinlein2016pit,heinlein2016pto, heinlein2016pos} of Trilinos~\cite{trilinos-website}.  We combine the implemented monolithic overlapping Schwarz preconditioner~\cite{heinlein2019monolithic,hochmuth2020parallel} and the multi-level implementation~\cite{roever2022multilevel:phd, heinlein2023multilevel} to get a parallel monolithic multi-level overlapping Schwarz preconditioner.  Our implementation is algebraic in the sense that we need only information about the parallel distribution of dofs of the velocity and pressure.  Everything else, e.g., the construction of the local subdomain matrices, the construction of restriction and interpolation operators, and the construction of the coarse problem, is handled by the FROSch package.

Furthermore, we implemented an interface between the FEAT software library~\cite{featflow-website}, which provides numerical assembly routines for finite element problems and a geometric multigrid solver package, see, e.g.~\cite{duennebacke2021space-time,ruda2022accelerator, paul2017preconditioning}, and the FROSch package.  We can assemble the system matrix and a right-hand side vector of a finite element problem in FEAT and use FROSch as a solver through the interface.  In the interface, the matrix and vector will be converted into a Trilinos Tpetra matrix and vector, which serve as in input for FROSch.  Note that FROSch can also be used as a coarse solver of the geometric multigrid solver of FEAT, which will not be discussed in this work.

Note that in our numerical experiments, we use Q2 elements for the velocity $u$ and P1-discontinuous elements for the pressure $p$.  Since the pressure is also discontinuous across the subdomains, we do not have an interface $\Gamma_{p}$ of the dofs of $p$.  Therefore, as coarse basis functions for the pressure we use constant functions over the whole subdomains, i.e., for each subdomain we have one coarse basis function of the pressure which is constant over the whole subdomain.

\section{Numerical Results}\label{pamm-25-koehler-sec:numerical-results}
For our numerical experiments, we consider a Poiseuille flow on the unit cube, to show the parallel performance of our implementation, and a Poiseuille flow and more advanced geometry, an extrusion die, which can be used in industrial applications.

For the Poiseuille flow on the unit cube, we use a parabolic inflow profile of type $u(x,y,z) = \left[cyz(1 - y)(1 - z), 0, 0\right]^{T}$, where $c$ is a constant.  The results of a weak scaling test are reported in Table~\ref{pamm-25-koehler-table:unit-cube-poiseuille-flow-2-lvl} and Table~\ref{pamm-25-koehler-table:unit-cube-poiseuille-flow-3-lvl}.  Note that the subdomains are cubes and therefore the decomposition is structured, see Figure~\ref{pamm-25-koehler-fig:unit-cube-decomposition} which is beneficial for domain decomposition methods due to their dependency of the ratio of the subdomain diameter to the element diameter.  As coarse space, we use a variant of the GDSW coarse space:  The GDSW$^{*}$ coarse space, see~\cite{hochmuth2020parallel}.  The difference is that here the vertices and adjacent edges are grouped to one interface entity.  This leads to a smaller coarse space and works here very well.

In Table~\ref{pamm-25-koehler-table:unit-cube-poiseuille-flow-2-lvl}, we report the results for the two-level solver.  The number of Krylov iteration increases from 34 to 37 when we go from 512 to 4\,096 subdomains, which is acceptable.  The symbolic set up time includes the construction of overlapping information for FROSch preconditioner and also the construction of matrix (without numerical values) which will be used by FROSch.  The symbolic set up time increases by nearly a factor of four, from 1.85s up to 7.78s.  Note that the most time consuming part here is the construction of the overlapping information, which is computed in an algebraic fashion.  This time can be decreased by providing the information which dof belongs to the overlap of which subdomain.  The numerical set up time includes the numerical initialization of the matrix with concrete values, the local factorization of the subdomains as well as the factorization of the coarse problem, and other data structures that are needed by FROSch.  This increases also by nearly a factor of 4.5, from 8.59 up to 38.2s, which is mostly due to factorization time of the direct solver for much larger coarse space, which has increased by nearly a factor of 9.
%The time for the solution, application of the Krylov method, 
The solve time (Krylov iteration)
increase from 0.69s to 1.84s, which can also be explained by the larger coarse space, and more time for the forward-backward solve due to more MPI communication
%, since we have factor 8 more subdomains.  
For 32\,768 subdomains, the factorization of the coarse problem (407\,837 dofs) fails.

\begin{vchtable}[t]
  \centering
  {
    \footnotesize
    \begin{tabular}{r r r  |  r r r r}
      \toprule
      \multicolumn{7}{c}{Monolithic \textbf{2-level} overlapping Schwarz Preconditioner: Poiseuille flow on the unit cube}                           \\[0.5ex]
        \#dof       &  \#subd. & {\#coarse dof}   & solve [in s]   & \multicolumn{2}{c}{setup} [in s]   & \#it.   \\
                    &          &    2nd-level     &                & symbolic       & numeric           &         \\\hline\rule{0pt}{3ex}
          954\,947  &     512  &    5\,573        & 0.69           &   1.85         &     8.59          &  34     \\
       7\,488\,643  &  4\,096  &   48\,781        & 1.84           &   7.78         &    38.2           &  37     \\
      59\,312\,387  & 32\,768  &  407\,837        &   $\infty$          &   $\infty$          &      $\infty$          &   $\infty$   \\
      \bottomrule
    \end{tabular}
  }
  \caption{\footnotesize 3$D$ Poiseuille flow on the unit cube solved with monolithic 2-level overlapping Schwarz:  discretization Q2-P1-discontinuous;  assembled in FEAT and solved with Trilinos/FROSch with GDSW$^{*}$ coarse space; weak scaling results; Krylov method:  FGMRES; stopping criterion: $\text{tol}_{\text{rel}}=10^{-8}$}\label{pamm-25-koehler-table:unit-cube-poiseuille-flow-2-lvl}
\end{vchtable}

\begin{vchfigure}
  \centering
  \begin{subfigure}[b]{0.48\textwidth}
    \centering
    \includegraphics[width=0.35\textwidth]{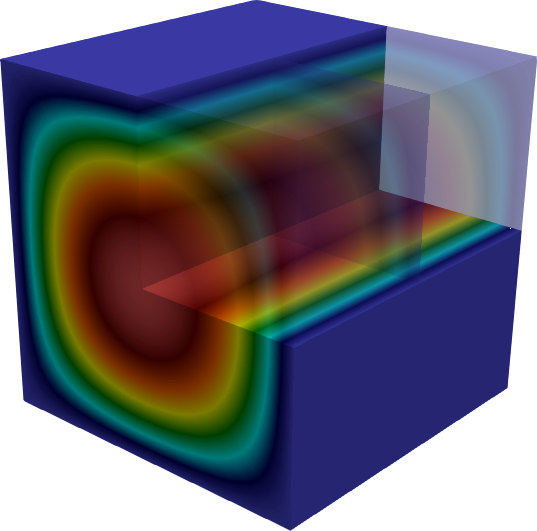}\par\vspace*{1ex}
    \includegraphics[scale=0.15]{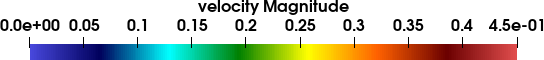}
    \caption{Magnitude of the velocity of a Poiseuille flow on the unit cube.  There is one inflow boundary and one outflow boundary.}\label{pamm-25-koehler-fig:unit-cube-solution}
  \end{subfigure}\hfill
  \begin{subfigure}[b]{0.48\textwidth}
    \includegraphics[width=0.35\textwidth]{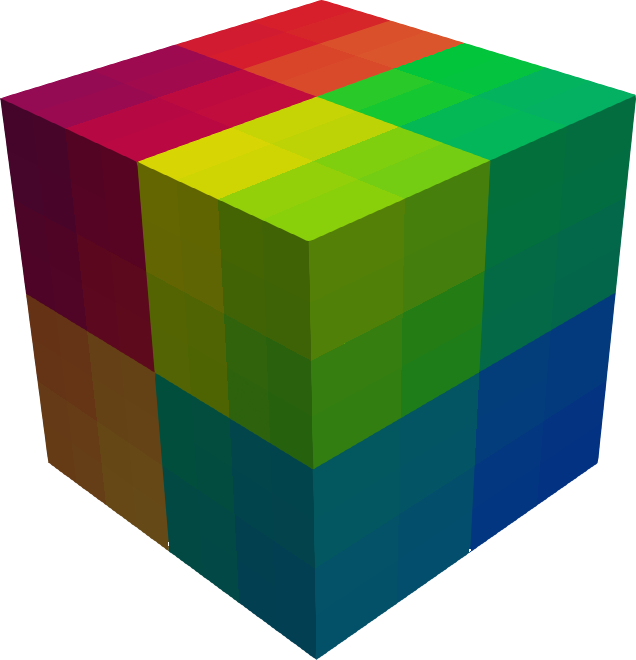} \hspace*{4ex}
    \includegraphics[width=0.35\textwidth]{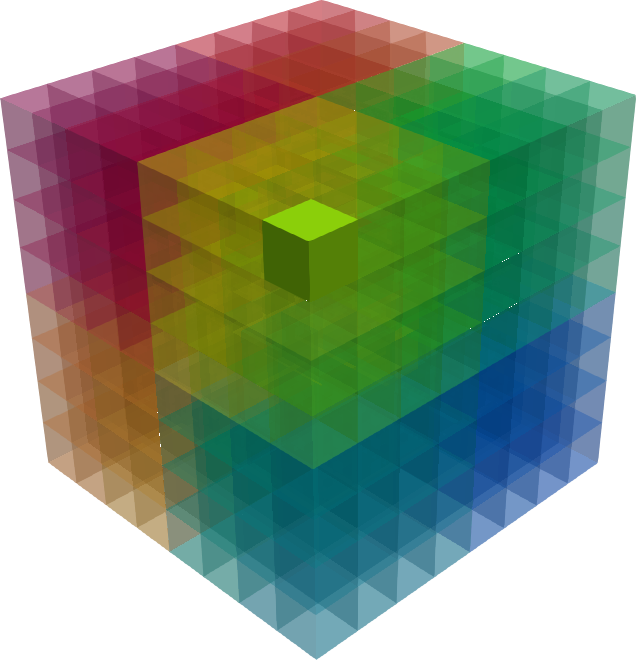}
    \caption{\emph{Left:}  All 512 subdomains are highlighted. \emph{Right:} Subdomain 220 is highlighted.}\label{pamm-25-koehler-fig:unit-cube-decomposition}
  \end{subfigure}
  \caption{Visualization of the solution of a Poiseuille flow on the unit cube and of the decomposition into 512 subdomains.}\label{pamm-25-koehler-fig:unit-cube}
\end{vchfigure}

In Table~\ref{pamm-25-koehler-table:unit-cube-poiseuille-flow-3-lvl}, we report the results for the monolithic three-level overlapping Schwarz solver.  Let us highlight that all three cases, 512, 4\,096, and 32\,768 subdomains, can be solved.  Furthermore, solution time and set up times for 4\,096 subdomains are faster than the timings for the two-level preconditioner.  The number coarse dofs on the second level denotes the size of $\mathcal{A}_{0}^{(2)}$, see~\eqref{eq:pamm-feat-frosch-gdsw-three-level}, which is the same as in Table~\ref{pamm-25-koehler-table:unit-cube-poiseuille-flow-2-lvl} and the coarse dofs on the third level denotes the size of $\mathcal{A}_{0}^{(3)}$.  With increasing number of subdomains, we also increase the number subregions (subdomains of the third level), which works the best here.  The number of Krylov iterations increase from 41 to 66, which is still acceptable, since we increase the number of subdomains by a factor 64.  One possibility for the increasing number of Krylov iterations might be that we are still not in the asymptotic range of the preconditioner.  The symbolic set up time increases much more, from 1.21s to 77.7s.  Note that, like in the two-level case, this is due to the computation of the overlapping information and construction of the matrix.  The difference for 512 and 4\,096 subdomains for the two-level and three-level is very small and due to different runs on different compute nodes.  The numerical set up time includes here, beside the necessary set up for the two-level data structures, also the partition of the coarse problem into subregions and the computation of the overlapping information for the subregions.  The time increases from 3.61s up to 8.87s, which is acceptable, since we 4\,213 coarse dofs on the third level, which is in same magnitude like the coarse problem as the coarse problem for the two-level approach for 512 subdomains, which has 5\,573 coarse dofs and needs 4.46s for the numerical set up time.

\begin{vchtable}[t]
  \centering
  {
    \footnotesize
    \begin{tabular}{r r r r r |  r r r r}
      \toprule
      \multicolumn{9}{c}{Monolithic \textbf{3-level} overlapping Schwarz Preconditioner: Poiseuille flow on the unit cube}                           \\[0.5ex]
        \#dof       &  \#subd. &  \#subreg. & \multicolumn{2}{c}{\#coarse dof}   & solve [in s]   & \multicolumn{2}{c}{setup} [in s]   & \#it.   \\
                    &          &            &    2nd-level   &  3rd-level        &                & symbolic       & numeric           &         \\\hline\rule{0pt}{3ex}
          954\,947  &     512  &       8    &    5\,573      &          65       &   0.65         &    1.21        &   3.73            &  41     \\
       7\,488\,643  &  4\,096  &      64    &   48\,781      &         901       &   1.26         &    7.42        &   3.61            &  52     \\
      59\,312\,387  & 32\,768  &     256    &  407\,837      &      4\,213       &   6.22         &   77.7         &   8.87            &  66     \\
      \bottomrule
    \end{tabular}
  }
  \caption{\footnotesize 3$D$ Poiseuille flow on the unit cube solved with monolithic 3-level overlapping Schwarz:  discretization Q2-P1-discontinuous;  assembled in FEAT and solved with Trilinos/FROsch with GDSW$^{*}$ coarse space; weak scaling results; Krylov method:  FGMRES; stopping criterion: $\text{tol}_{\text{rel}}=10^{-8}$}\label{pamm-25-koehler-table:unit-cube-poiseuille-flow-3-lvl}
\end{vchtable}

In Table~\ref{pamm-25-koehler-table:unit-cube-poiseuille-flow-2-lvl} and Table~\ref{pamm-25-koehler-table:unit-cube-poiseuille-flow-3-lvl}, we report results for a weak scaling test for a Poiseuille flow example with two inflow boundaries and one outflow boundary, see Figure~\ref{pamm-25-koehler-fig:gendie}.  Note that the geometry is quite more complex than the unit cube and decomposition is not structured, see Figure~\ref{pamm-25-koehler-fig:gendie-decomposition}.

\begin{vchtable}[t]
  \centering
  {
    \footnotesize
    \begin{tabular}{r r r  |  r r r r}
      \toprule
      \multicolumn{7}{c}{Monolithic \textbf{2-level} overlapping Schwarz Preconditioner: Poiseuille flow on an extrusion die}                           \\[0.5ex]
        \#dof       &  \#subd. & {\#coarse dof}   & solve [in s]   & \multicolumn{2}{c}{setup} [in s]   & \#it.   \\
                    &          &    2nd-level     &                & symbolic       & numeric           &         \\\hline\rule{0pt}{3ex}
       2\,358\,567  &     512  &     8\,408       &  2.36          &   2.06         &    11.6           &  67     \\
      17\,603\,097  &  4\,096  &    80\,977       &  9.63          &  16.8          &    84.3           &  61     \\
      \bottomrule
    \end{tabular}
  }
  \caption{\footnotesize 3$D$ Poiseuille flow on an extrusion die solved with monolithic 2-level overlapping Schwarz :  discretization Q2-P1-discontinuous;  assembled in FEAT and solved with Trilinos/FROsch with GDSW$^{*}$ coarse space; weak scaling results; Krylov method:  FGMRES; stopping criterion: $\text{tol}_{\text{rel}}=10^{-10}$}\label{pamm-25-koehler-table:gendie-poiseuille-flow-2-lvl}
\end{vchtable}

\begin{vchtable}[t]
  \centering
  {
    \footnotesize
    \begin{tabular}{r r r r r |  r r r r}
      \toprule
      \multicolumn{9}{c}{Monolithic \textbf{3-level} overlapping Schwarz Preconditioner: Poiseuille flow on an extrusion die}                           \\[0.5ex]
        \#dof       &  \#subd. &  \#subreg. & \multicolumn{2}{c}{\#coarse dof}   & solve [in s]   & \multicolumn{2}{c}{setup} [in s]   & \#it.   \\
                    &          &            &    2nd-level   &  3rd-level        &                & symbolic       & numeric           &         \\\hline\rule{0pt}{3ex}
       2\,358\,567  &     512  &       8    &    8\,408      &         161       &   2.59         &    2.08        &   8.83            &  82     \\
      17\,603\,097  &  4\,096  &      64    &   80\,977      &      1\,102       &   5.12         &   16.7         &  14.0             &  95     \\
      \bottomrule
    \end{tabular}
  }
  \caption{\footnotesize 3$D$ Poiseuille flow on an extrusion die solved with monolithic 3-level overlapping Schwarz:  discretization Q2-P1-discontinuous;  assembled in FEAT and solved with Trilinos/FROsch with GDSW$^{*}$ coarse space; weak scaling results; Krylov method:  FGMRES; stopping criterion: $\text{tol}_{\text{rel}}=10^{-10}$}\label{pamm-25-koehler-table:gendie-poiseuille-flow-3-lvl}
\end{vchtable}

In Table~\ref{pamm-25-koehler-table:unit-cube-poiseuille-flow-2-lvl}, we report results for the monolithic two-level overlapping Schwarz preconditioner.  The same technical details as before apply.  The symbolic set up increase when we increase the number of subdomains from 512 to 4\,096 subdomains due to greater cost for the computation of the overlapping information and the construction of the matrix.  The numerical set up time increases from 11.6s up to 84.3s due to the larger coarse space.  The solve time increases also mostly due to the larger coarse space.  The number of Krylov iterations decreases which is an issue with the relative stopping tolerance of the Krylov method and the increasing initial residual.

In Table~\ref{pamm-25-koehler-table:unit-cube-poiseuille-flow-3-lvl}, we report results for the monolithic three-level overlapping Schwarz preconditioner.  As in the comparison of the two- and three-level preconditioner for the unit cube, the set up time decreases for the three-level preconditioner.  Indeed, the numerical set up time decreases about a factor of 6 and the symbolic set up stays the same.  The time for the solution decreases nearly about a factor of 2.

\begin{vchfigure}
  \centering
  \includegraphics[width=0.17\textwidth]{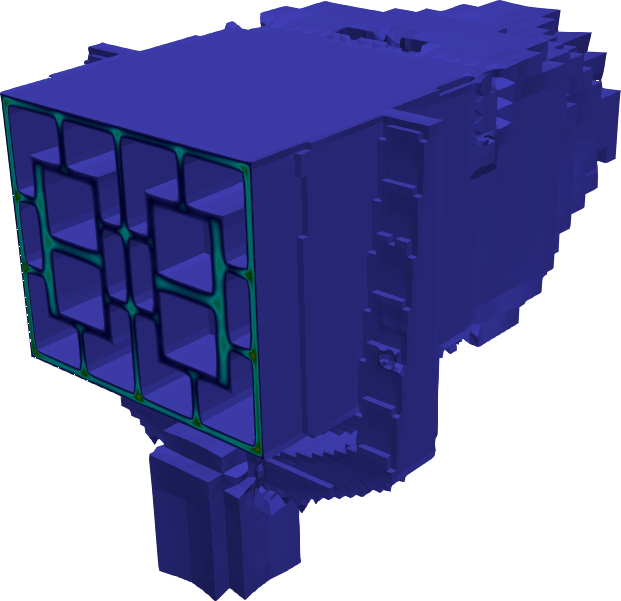} \hspace*{3ex}
  \includegraphics[width=0.13\textwidth]{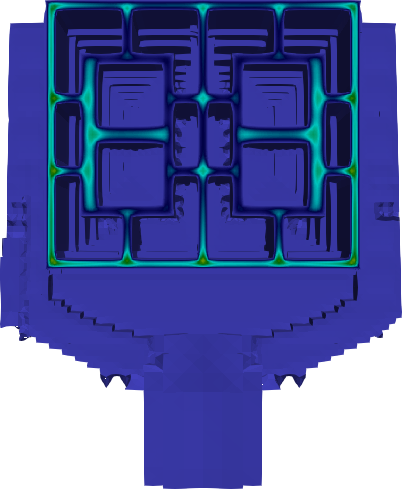} \hspace*{3ex} %\vspace*{1ex}
  \includegraphics[width=0.15\textwidth]{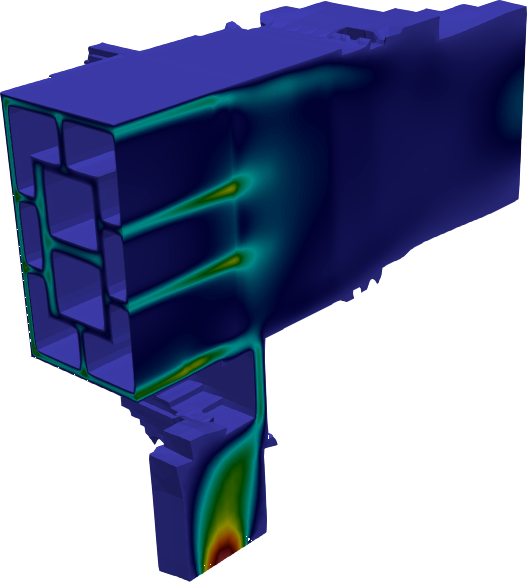} \hspace*{3ex}
  \includegraphics[width=0.17\textwidth]{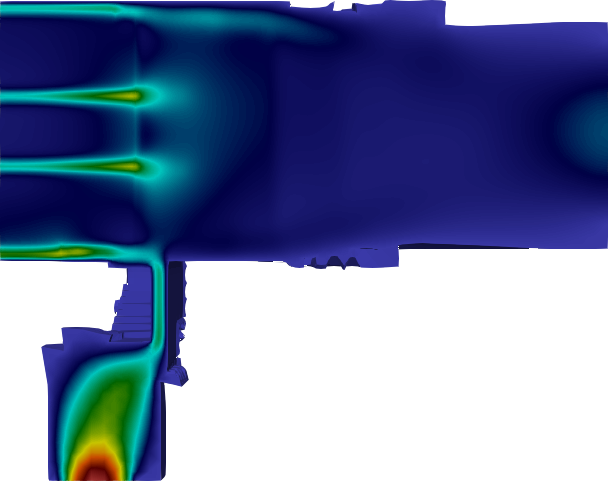} \vspace*{2ex}
  \par
  \includegraphics[scale=0.2]{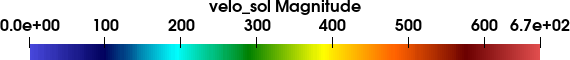}
  \caption{Visualization of the magnitude of the velocity of a Poiseuille flow on an extrusion die.  There are two inflow boundaries and one outflow boundary, see first image from the right: left and bottom boundary are the inflow boundaries and right boundary is the outflow boundary.}\label{pamm-25-koehler-fig:gendie}
\end{vchfigure}

\begin{vchfigure}
  \centering
  \begin{subfigure}[t]{0.49\textwidth}
    \centering
    \includegraphics[width=0.35\textwidth]{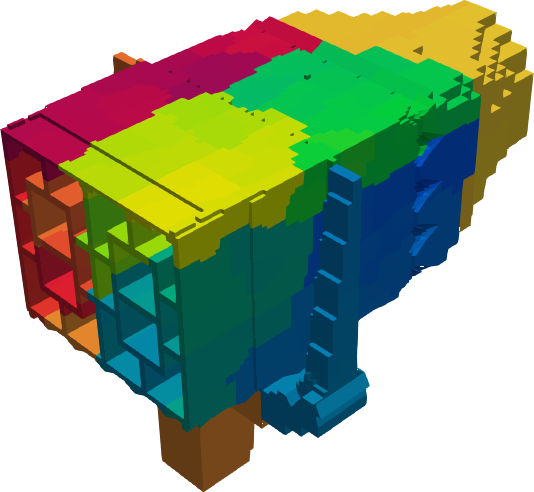} \hspace*{4ex}
    \includegraphics[width=0.22\textwidth]{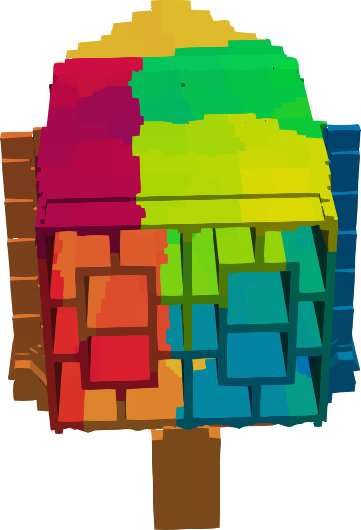} % \hspace*{4ex}
    \caption{All 512 subdomains are highlighted.}
  \end{subfigure}
  \begin{subfigure}[t]{0.49\textwidth}
    \includegraphics[width=0.35\textwidth]{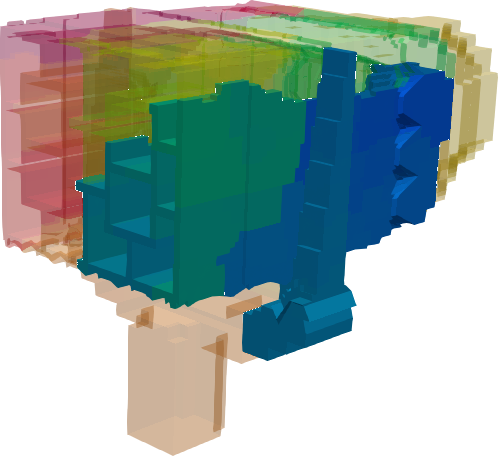} \hspace*{4ex}
    \includegraphics[width=0.35\textwidth]{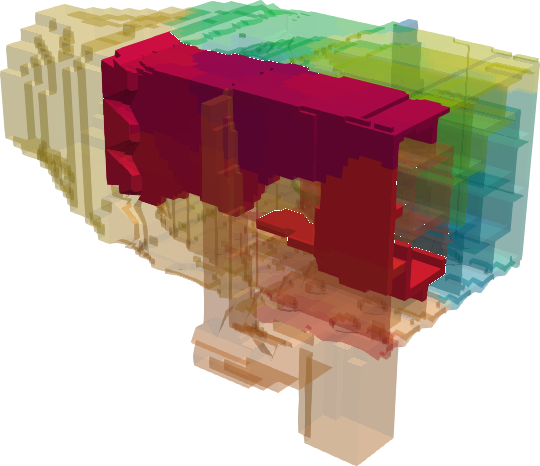}
    \caption{\emph{Left:}  Subdomains 1--129 are highlighted. \emph{Right:} Subdomains 433--512 are highlighted.}
  \end{subfigure}
  \caption{Visualization of the decomposition into 512 subdomains of the extrusion die.}\label{pamm-25-koehler-fig:gendie-decomposition}
\end{vchfigure}

With the three-level preconditioner, we can solve much larger problems, see Table~\ref{pamm-25-koehler-table:unit-cube-poiseuille-flow-2-lvl} and Table~\ref{pamm-25-koehler-table:unit-cube-poiseuille-flow-3-lvl} for the case of 32\,768 subdomains, and obtain a performance gain about of factor 6 in the set up time, see Table~\ref{pamm-25-koehler-table:gendie-poiseuille-flow-2-lvl} and Table~\ref{pamm-25-koehler-table:gendie-poiseuille-flow-3-lvl}.

\section{Conclusion}
We implemented successfully a monolithic three-level overlapping Schwarz preconditioner within the Trilinos/FROSch software package and tested the scalability with two examples:  A Stokes flow problem on the unit cube and on a more complex geometry of an extrusion die.  Our numerical experiments show the benefit of the three-level preconditioner compared to the two-level preconditioner, when the coarse space grows larger and larger.

\begin{acknowledgement}
  This work is funded by the German Federal Ministry of Education and Research under grant no. 16ME0708.

  Funded by the European Union - NextGenerationEU.

  The authors gratefully acknowledge the scientific support and HPC resources provided by the Erlangen National High Performance Computing Center (NHR@FAU) of the Friedrich-Alexander-Universität Erlangen-Nürnberg (FAU) under the NHR project k107ce. NHR funding is provided by federal and Bavarian state authorities. NHR@FAU hardware is partially funded by the German Research Foundation (DFG) – 440719683.

  We would like to thank Peter Zajac and Maximilian Esser from TU Dortmund for their part of the implementation of the FEAT and FROSch interface and the implementation of the numerical examples.
\end{acknowledgement}

\vspace{\baselineskip}
%% The style of the following references should be used in all documents.

\bibliographystyle{pamm}

\providecommand{\WileyBibTextsc}{}
\let\textsc\WileyBibTextsc
\providecommand{\othercit}{}
\providecommand{\jr}[1]{#1}
\providecommand{\etal}{~et~al.}

\end{document}